\documentclass [11pt,oneside]{article}

\reversemarginpar \textwidth=14cm \textheight=19cm

\usepackage{amssymb} \usepackage{amsfonts} \usepackage{amsmath}
\usepackage{amsthm} \usepackage{epsfig}

\newcommand{\mettifig}[1]{\epsfig{file=#1}}

\newtheorem{lemma}{Lemma}[section] 

\newtheorem{rem}[lemma]{Remark} 

\newtheorem{cor}[lemma]{Corollary}

\newcommand{\matR} {\ensuremath {\mathbb{R}}}
\newcommand{\matQ} {\ensuremath {\mathbb{Q}}}
\newcommand{\matZ} {\ensuremath {\mathbb{Z}}}

\newcommand{\matH} {\ensuremath {\mathbb{H}}}

\newcommand{\calK} {\ensuremath {\mathcal{K}}}

\newcommand{\calM} {\ensuremath {\mathcal{M}}}
\newcommand{\calC} {\ensuremath {\mathcal{C}}}

\newcommand{\calP} {\ensuremath {\mathcal{P}}}
\newcommand{\calT} {\ensuremath {\mathcal{T}}}

\newcommand{\calH}{\ensuremath {\mathcal{H}}}

\newcommand{\nota} [1] {\caption{\footnotesize{#1}}}
\newcommand{\matr} [4] {\left(\begin{array}{cc} #1 & #2 \\ #3 & #4 \\ \end{array} \right)}

\newfont{\Got}{eufm10 scaled 1200}

\font\titsc=cmcsc10 scaled 1200

\newcommand{\diag}{\SelectTips{cm}{} \xymatrix@1}

\newcommand{\timtil}{\begin{picture}(12,12)
\put(2,0){$\times$}\put(2,4.5){$\sim$}\end{picture}}

\author{Roberto \titsc{Frigerio} \and Bruno \titsc{Martelli}
\and Carlo \titsc{Petronio}}

\title{Small hyperbolic 3-manifolds with geodesic boundary}

\begin{document}

\maketitle

\begin{abstract}
\noindent
We classify the orientable finite-volume hyperbolic 3-manifolds having 
non-empty compact totally geodesic boundary and admitting
an ideal triangulation with at most four tetrahedra.
We also compute the volume of all such manifolds, we describe
their canonical Kojima decomposition, and we discuss manifolds 
having cusps. 

The manifolds built from one or two tetrahedra were previously known. There are 151 different manifolds
built from three tetrahedra, realizing 18 different volumes. Their Kojima decomposition always consists
of tetrahedra (but occasionally requires four of them). And there is a single cusped manifold, that we
can show to be a knot complement in a genus-2 handlebody. Concerning manifolds built from four
tetrahedra, we show that there are 5033 different ones, with 262 different volumes.  The Kojima
decomposition consists either of tetrahedra (as many as eight of them in some cases), or of two pyramids,
or of a single octahedron. There are 30 manifolds having a single cusp, and one having two cusps.

Our results were obtained with the aid of a computer.
The complete list of manifolds (in SnapPea format) 
and full details on their invariants are
available on the world wide web.

  \vspace{4pt}

\end{abstract}

\noindent
This paper is devoted to the class of all orientable finite-volume hyperbolic 
3-manifolds having non-empty compact totally geodesic boundary 
and admitting a minimal ideal triangulation with either three 
or four but no fewer tetrahedra. We describe the theoretical background
and experimental results of a computer program that has enabled us to
classify all such manifolds. (The case of manifolds obtained from
two tetrahedra was previously dealt with in~\cite{Fuj}).
We also provide an overall discussion of the most important features
of all these manifolds, namely of:
\begin{itemize}
\item their volumes;
\item the shape of their canonical Kojima decomposition;
\item the presence of cusps.
\end{itemize}
These geometric invariants have all been determined by 
our computer program. The complete list of manifolds in SnapPea format and 
detailed information on the invariants is available from~\cite{www}.

\section{Preliminaries and statements}\label{statements:section}
We consider in this paper the class $\calH$ of orientable 3-manifolds $M$
having compact non-empty boundary $\partial M$ and admitting
a complete finite-volume hyperbolic metric with respect to which
$\partial M$ is totally geodesic. It is a well-known fact~\cite{kojima} that
such an $M$ is the union of a compact portion and some cusps
based on tori, so it has a natural
compactification obtained by adding some tori.
The elements of $\calH$ are regarded up to homeomorphism, or
equivalently isometry (by Mostow's rigidity).

\paragraph{Candidate hyperbolic manifolds}
Let us now introduce the class $\widetilde\calH$ of 3-manifolds $M$ such that:
\begin{itemize}
\item $M$ is orientable, compact, boundary-irreducible and acylindrical;
\item $\partial M$ consists of some tori (possibly none of them) and
at least one surface of negative Euler characteristic.
\end{itemize}
The basic theory of hyperbolic manifolds implies that,
up to identifying a manifold with its natural compactification,
the inclusion $\calH\subset\widetilde{\calH}$ holds.
We note that, by Thurston's hyperbolization, an element of $\widetilde\calH$
actually lies
in $\calH$ if and only if it is atoroidal.
However we do not require atoroidality in the definition of $\widetilde\calH$,
for a reason that will be mentioned later in this section and explained in detail
in Section~\ref{spines:section}.

Let $\Delta$ denote the standard tetrahedron, and let $\Delta\!^*$ be
$\Delta$ minus open stars of its vertices. 
Let $M$ be a compact $3$-manifold with $\partial M\neq\emptyset$.
An \emph{ideal triangulation} of $M$
is a realization of $M$ as a gluing of a finite
number of copies of $\Delta\!^*$, induced by a simplicial face-pairing  
of the corresponding $\Delta$'s. We denote by $\calC_n$ the class of
all orientable manifolds admitting an ideal triangulation with
$n$, but no fewer, tetrahedra, and we set:
$$\calH_n=\calH\cap\calC_n,\qquad \widetilde{\calH}_n=\widetilde{\calH}\cap\calC_n.$$
We can now quickly explain why we did not include atoroidality in
the definition of $\widetilde{\calH}$. The point is that there is
a general notion~\cite{Mat} of \emph{complexity} $c(M)$ for a 
compact 3-manifold $M$, and
$c(M)$ coincides with the minimal number of tetrahedra in an ideal
triangulation precisely when $M$ is boundary-irreducible and 
acylindrical. This property makes it feasible to
enumerate the elements of $\widetilde{\calH}_n$. 

To summarize our definitions, \emph{we can interpret $\calH_n$ as the
set of $3$-manifolds which have complexity $n$ and are hyperbolic
with non-empty compact geodesic boundary, while $\widetilde{\calH}_n$
is the set of complexity-$n$ manifolds which are only ``candidate
hyperbolic.''} 

\paragraph{Enumeration strategy}
The general strategy of our classification result
is then as follows:
\begin{itemize}
\item We employ the technology of standard spines~\cite{Mat} (and more 
particularly o-graphs~\cite{BePe}), together with certain 
\emph{minimality tests}
(see Section~\ref{spines:section} below), to produce for $n=3,4$ a list of 
triangulations
with $n$ tetrahedra such that every element of $\widetilde{\calH}_n$
is represented by some triangulation in the list.
Note that the same element of $\widetilde{\calH}_n$ is represented
by several distinct triangulations. Moreover, there could \emph{a priori}
be in the list triangulations
representing manifolds of complexity lower than $n$, but
the result of the classification itself actually shows that
our minimality tests are sophisticated enough to ensure
this does not happen; 
\item We write and solve the hyperbolicity equations (see~\cite{FriPe}
and Section~\ref{hyperbolicity:section} below)
for all the triangulations, finding solutions in the vast majority
of cases (all of them for $n=3$);
\item We compute the tilts (see~\cite{FriPe,Ush} 
and Section~\ref{hyperbolicity:section} again)
of each of the geometric 
triangulations thus found,
whence determining whether the triangulation (or maybe a partial
assembling of the tetrahedra of the triangulation) gives
Kojima's canonical decomposition; when it does not, we 
modify the triangulation according to the strategy described in~\cite{FriPe},
eventually finding the canonical decomposition in all cases;
\item We compare the canonical decompositions to each other, thus finding
precisely which pairs of triangulations in the list represent identical manifolds;
we then build a list of distinct hyperbolic manifolds, which coincides with
$\calH_n$ because of the next point;
\item We prove that when the hyperbolicity equations have no solution then indeed
the manifold is not a member of $\calH_n$, because it contains an 
incompressible torus (this is shown in Section~\ref{spines:section}).
\end{itemize}
Even if the next point is not really part of the classification strategy, 
we single it out as an important one:
\begin{itemize}
\item We compute the volume of all the elements of $\calH_n$ using
the geometric triangulations already found and the formulae from~\cite{ushi-vol}.
\end{itemize}

\paragraph{One-edged triangulations}
Before turning to the description of our discoveries, we must mention another
point. Let us denote by $\Sigma_g$ the orientable surface of genus $g$, and
by $\calK(M)$ the blocks of the canonical Kojima decomposition of $M\in\calH$.
We have introduced in~\cite{FriMaPe1} the class $\calM_n$ of orientable manifolds
having an ideal triangulation with $n$ tetrahedra and a single edge, and we have shown
that for $n\geq2$ and $M\in\calM_n$:
\begin{itemize}
\item $M$ is hyperbolic with geodesic boundary $\Sigma_n$; 
\item $M$ has a unique ideal triangulation with $n$ tetrahedra, which coincides
with $\calK(M)$; moreover $c(M)=n$ and $\calM_n=\{M\in\widetilde\calH_n:\ \partial M=\Sigma_n\}$;
\item the volume of $M$ depends only on $n$ and can be computed explicitly.
\end{itemize}
These facts imply in particular that $\calM_n$ is contained in $\calH_n$. 

\paragraph{Results}
We can now state our main results, recalling first~\cite{Fuj} 
that $\calH_1=\emptyset$ and $\calH_2=\calM_2$ has eight elements,
and pointing out that
all the values of volumes in our statements are approximate, not
exact ones.  More accurate approximations are available on the web~\cite{www}.
We also emphasize that our results indeed have an experimental nature,
but we have checked by hand a number of cases and always found perfect
agreement with the results found by the computer.

\paragraph{Results in complexity 3} We have discovered that:
\emph{\begin{itemize}
\item $\calH_3$ coincides with $\widetilde{\calH}_3$ and has $151$ elements; 
\item $\calM_3$ consists of $74$ elements of volume $10.428602$;
\item all the $77$ elements of $\calH_3\setminus\calM_3$ have 
boundary $\Sigma_2$, and one of them also has one cusp.
\end{itemize}
\noindent Moreover the elements $M$ of $\calH_3\setminus\calM_3$ split as follows:
\begin{itemize} 
\item $73$ compact $M$'s with $\calK(M)$ consisting of three tetrahedra; ${\rm vol}(M)$
attains on them $15$ different values,
ranges from $7.107592$ to $8.513926$, and has maximal multiplicity nine,
with distribution according to 
number of manifolds as shown in Table~\ref{c3:table} (see the Appendix);
\item three compact $M$'s with $\calK(M)$ consisting of 
four tetrahedra; they all have the same volume $7.758268$;
\item one non-compact $M$; it has a single toric 
cusp, $\calK(M)$ consists of three tetrahedra, and ${\rm vol}(M)=7.797637$.
\end{itemize}}

The cusped element of $\calH_3$ turns out to be a very interesting manifold.
In~\cite{FriMaPe3} we have analyzed all the Dehn fillings of its toric 
cusp,
improving previously known bounds on the distance between non-hyperbolic
fillings. In particular, we have shown that there are fillings giving the
genus-$2$ handlebody, so the manifold in question is a knot complement, as shown
in Fig.~\ref{cusp_is_complement:fig}.
\begin{figure}
\begin{center}
\mettifig{genus2,width=7cm}
\nota{The cusped manifold having complexity three and non-empty
boundary is the complement of a knot in the
genus-two handlebody.} \label{cusp_is_complement:fig}
\end{center}
\end{figure}

\paragraph{Results in complexity 4} We have discovered that:
\emph{\begin{itemize}
\item $\calH_4$ has $5033$ elements, and $\widetilde{\calH}_4$ and has $6$ more;
\item $5002$ elements of $\calH_4$ are compact; more precisely:
\begin{itemize}
\item $2340$ have boundary $\Sigma_4$ (\emph{i.e.}~they belong to $\calM_4$);
\item $2034$ have boundary $\Sigma_3$;
\item $628$ have boundary $\Sigma_2$;
\end{itemize}
\item $31$ elements of $\calH_4$ have cusps; more precisely:
\begin{itemize}
\item $12$ have one cusp and boundary $\Sigma_3$;
\item $18$ have one cusp and boundary $\Sigma_2$;
\item one has two cusps and boundary $\Sigma_2$.
\end{itemize}
\end{itemize}}

More detailed information about the volume and the shape of the
canonical Kojima decomposition of these manifolds is described 
in Tables~\ref{c4:cpt:summary:table} and~\ref{c4:cusp:summary:table}.
In these tables each box corresponds to the 
manifolds $M$ having a prescribed boundary and type of $\calK(M)$.
The first information we provide (in boldface)
within the box is the number of distinct such $M$'s.
When all the $M$'s in the box have the same volume, we 
indicate its value. 
Otherwise we indicate the minimum, the maximum,
the number of different values, and the maximal multiplicity of the values
of the volume function, and we refer to one of the 
tables in the Appendix where more accurate 
information can be found. We emphasize here that, just as above,
$\calK(M)$ only describes the \emph{blocks} of the Kojima decomposition, not the 
combinatorics of the gluing.

\begin{table}
\begin{center}
\begin{tabular}{l||c|c|c}
& $\Sigma_4$ & $\Sigma_3$ & $\Sigma_2$ \\ \hline\hline
4 tetra & 
	{\bf 2340} & {\bf 1936} & {\bf 555} \\
	& 
	${\rm vol}=14.238170$ & $\min({\rm vol})=11.113262$ & $\min({\rm vol})=7.378628$ \\	
	& & $\max({\rm vol})=12.903981$ & $\max({\rm vol})=10.292422$ \\
	& & ${\rm values}=59$ & ${\rm values}=169$ \\
	& & ${\rm max\ mult}=138$ & ${\rm max\ mult}=27$ \\
	& & (Tables~\ref{c4:bd3:K4:first:table} and~\ref{c4:bd3:K4:second:table}) & 
		(Tables~\ref{c4:bd2:K4:first:table} and~\ref{c4:bd2:K4:second:table}) \\ \hline
5 tetra & & {\bf 42} & {\bf 41} \\
	& & ${\rm vol}=11.796442$ & $\min({\rm vol})=8.511458$ \\ 
	& & & $\max({\rm vol})=9.719900$\\
	& & & ${\rm values}=16$\\ 
	& & & ${\rm max\ mult}=6$\\
	& & & (Table~\ref{c4:bd2:K5:table}) \\ \hline
6 tetra & & & {\bf 3} \\ 
	& & & ${\rm vol}=8.297977$  \\ \hline
8 tetra & & & {\bf 3} \\ 
	& & & ${\rm vol}=8.572927$  \\ \hline
1 octa & & {\bf 56} & {\bf 14} \\
(regular) & & ${\rm vol}=11.448776$ & ${\rm vol}=9.415842$ \\ \hline
1 octa & & & {\bf 8} \\
(non-reg) & & & ${\rm vol}=8.739252$ \\ \hline
2 square & & & {\bf 4} \\
pyramids & & & ${\rm vol}=9.044841$ 
\end{tabular}
\nota{Number of compact elements of $\calH_4$, subdivided 
according to the boundary (columns) and shape of the canonical 
Kojima decomposition (rows); `tetra' and `octa' mean 
`tetrahedron' and `octahedron'
respectively, and `square pyramid' means `pyramid with
square basis.'}\label{c4:cpt:summary:table}
\end{center}\end{table}

\begin{table}
\begin{center}
\begin{tabular}{c||c|c|c}
& 1\ cusp, $\Sigma_3$
& 1\ cusp, $\Sigma_2$
& 2\ cusps, $\Sigma_2$
\\ \hline\hline
4 tetra & {\bf 12} & {\bf 16} & {\bf 1} \\
	& ${\rm vol}= 11.812681$ & $\min({\rm vol})=8.446655$ & ${\rm vol}=9.134475$ \\ 
	& & $\max({\rm vol})=9.774939$ & \\
	& & ${\rm values}=8$ & \\ 
	& & ${\rm max\ mult}=3$ & \\
	& & (Table~\ref{c4:cusp:table}) & \\ \hline
2 square & & {\bf 2} & \\
pyramids & & ${\rm vol}=8.681738$ &
\end{tabular}
\nota{Number of cusped elements of $\calH_4$, subdivided 
according to cusps and boundary (columns), and the shape of the canonical 
Kojima decomposition (rows).}\label{c4:cusp:summary:table}
\end{center}\end{table}

In addition to what is described in the tables, we have the
following extra information on the geometric shape of $\calK(M)$ when
it is given by an octahedron:
\begin{itemize}
\item the group of $56$ manifolds in Table~\ref{c4:cpt:summary:table}
is built from an octahedron with all 
dihedral angles equal to $\pi/6$; 
\item the group of $14$ manifolds in Table~\ref{c4:cpt:summary:table}
is built from an octahedron with all dihedral angles equal to $\pi/3$; 
\item the group of $8$ manifolds in Table~\ref{c4:cpt:summary:table}
is built from an octahedron with three dihedral angles
$2\pi/3$ along a triple of pairwise disjoint edges, and two more
complicated angles (one repeated 3 times, one 6 times).
\end{itemize}

A careful analysis of the values of volumes found leads to
the following consequences:

\begin{cor}
For $n=3,4$, the maximum of the volume on $\calH_n$ 
is attained at the elements of $\calM_n$.
\end{cor}

\begin{rem}
\emph{With the only exceptions discussed below in Remarks~\ref{volume:rem:1} and~\ref{volume:rem:2}, 
if two manifolds in $\calH_3\cup\calH_4$ have the same volume then they also have
the same complexity, boundary, and number of cusps.
Moreover, they typically also have the same 
geometric shape of the blocks of the 
Kojima decomposition (but of course not the same combinatorics of gluings).}
\end{rem}

\begin{rem}
\emph{There are 280 distinct 
values of volume we have found in our census, and the vast majority
of them correspond to more than one manifold.  As a matter of fact,
only 25 values are attained just once: 22 are in
Tables~\ref{c4:bd2:K4:first:table}
and~\ref{c4:bd2:K4:second:table},
two in Table~\ref{c4:cusp:table}, and one is the volume
of the cusped element of $\calH_3$.}
\end{rem}

\begin{rem}\label{volume:rem:1}
\emph{As stated above, there are three elements of 
$\calH_3$ with canonical decomposition made of four tetrahedra. The 
set of geometric shapes of these four tetrahedra is
actually the same in all three cases, and it turns out that the same
tetrahedra can also be glued to give five different elements of 
$\calH_4$. This gives the only example we have of elements 
$\calH_3$ having the same volume as elements of $\calH_4$.
The volume in question is 7.758268.}
\end{rem}

\begin{rem}\label{volume:rem:2}
\emph{The double-cusped manifold in $\calH_4$ has the same volume
$9.134475$ as two of the single-cusped ones (see
Table~\ref{c4:cusp:table}),
and it is probably worth mentioning a heuristic explanation for this fact.
Recall first that an ideal 
triangulation of a manifold induces a triangulation 
of the basis of the cusps. For $28$ of the single-cusped manifolds in $\calH_4$ this 
triangulation involves two triangles, but for two of them it involves four,
just as it does with the double-cusped manifold (both tori contain
two triangles).
In addition, the geometric shapes of the four triangles are the same in
all three cases. In other words, one sees here that four Euclidean triangles
can be used to build either two ``small'' Euclidean tori or
a single ``big'' Euclidean torus (in two different ways). 
So, in some sense, the three manifolds in question have the same
``total cuspidal geometry'' (even if two manifolds have one cusp and one has two).
This phenomenon already occurs in the case of manifolds 
without boundary~\cite{SnapPea}, and also in this case 
leads to equality of volumes. In the present case equality 
is also explained by the fact that 
the three manifolds in question have Kojima decomposition with the
same geometric shape of the blocks. In fact, each of them is the
gluing of four 
isometric partially truncated tetrahedra with three dihedral angles $\pi/3$
and three $\pi/6$. }
\end{rem}

The next information may also be of some interest:

\begin{rem}\label{toroidal:rem}
\emph{We will show below that 
the six manifolds in $\widetilde\calH_4\setminus\calH_4$ split along an
incompressible torus
into two blocks, one homeomorphic to 
the twisted interval bundle over the Klein bottle and
the other one to the cusped manifold that belongs to $\calH_3$.
These blocks give the JSJ decomposition of the manifolds involved.
We will also show that the manifolds are indeed distinct by analyzing the gluing
matrix of the JSJ decomposition.}
\end{rem}

\begin{rem}
\emph{As an ingredient of our arguments, we have completely classified
the combinatorially inequivalent ways of building an orientable manifold
by gluing together in pairs the faces of an octahedron. This topic was 
already mentioned in~\cite{bibbia} as an example of how difficult
classifying 3-manifolds could be (note that there are as many as $8505$
gluings to be compared for combinatorial equivalence).
For instance, the group of $56$ manifolds that
appears in Table~\ref{c4:cpt:summary:table} 
arises from the gluings of the octahedron such that
all the edges get glued together. The groups of 14 and 8 arise similarly,
requiring two edges and restrictions on their valence.}
\end{rem}

\begin{rem}
\emph{We have never included information about 
homology, because this invariant typically gives a much coarser 
information than the geometric invariants we have computed
(only 14 different homology groups arise for our 5184 manifolds).
We note however that it occasionally happens that two
manifolds having the same complexity, boundary, volume, and 
geometric blocks of the
canonical decomposition have different homology. The homology groups
we have found are $\matZ^2\oplus\matZ/_n$ for $n=1,\ldots,8$,
$\matZ^3\oplus\matZ/_n$ for $n=1,2,3,5$, $\matZ^4$, and 
$\matZ^2\oplus\matZ/_2\oplus\matZ/_2$.}
\end{rem}

\begin{rem}\label{compact:equations:rem}
\emph{Even if we have not yet introduced the hyperbolicity equations
that we use to find the geometric structures, we point out a remarkable
experimental discovery. The equations to be used in the cusped case are
qualitatively different (and a lot more complicated) than those to be used in
the compact case.  However, for all the $32$ cusped manifolds of the census,
the hyperbolic structure was first found as a limit of approximate solutions
of the compact equations.}
\end{rem}

\begin{rem}
\emph{For each $M$ in $\calH_3\cup\calH_4$, each of the (often multiple)
minimal triangulations of $M$ has been found to be geometric,
\emph{i.e.}~the corresponding set of hyperbolicity equations has been
proved to have a genuine solution. This strongly supports the conjecture
that ``minimal implies geometric,'' that one could already guess from
the cusped case~\cite{SnapPea}.}
\end{rem}

\begin{rem}
\emph{For each $M$ in $\calH_3\cup\calH_4$, the Kojima decomposition has been obtained
by merging some tetrahedra of a geometric triangulation of $M$. It follows
that the Kojima decomposition of every manifold in $\calH_3\cup\calH_4$ admits
a subdivision into tetrahedra.}
\end{rem}

\section{Spines and the enumeration method}\label{spines:section}
If $M$ is a compact orientable $3$-manifold,
let $t(M)$ be the minimal number of tetrahedra 
in an ideal triangulation of either 
$M$, when $\partial M\ne\emptyset$, or
$M$ minus any number of balls, when $M$ is closed.
The function $t$ thus defined
has only one nice property: it is finite-to-one.
In~\cite{Mat} Matveev has introduced another function $c$, 
which he called \emph{complexity}, having
many remarkable properties not satisfied by $t$. For instance, $c$ is
additive on connected sums,
and it does not increase when cutting along an incompressible surface.
Moreover it was proved in~\cite{Mat, Mat-Exp} that 
$c$ equals $t$ on the most interesting
3-manifolds, namely $c(M)=t(M)$ when $M$ is $\partial$-irreducible
and acylindrical, and $c(M)<t(M)$ otherwise. Therefore, if $\chi(M)<0$, we have
$c(M)=t(M)$ if and only if $M\in\widetilde\calH$.

\begin{figure}
\begin{center}
\mettifig{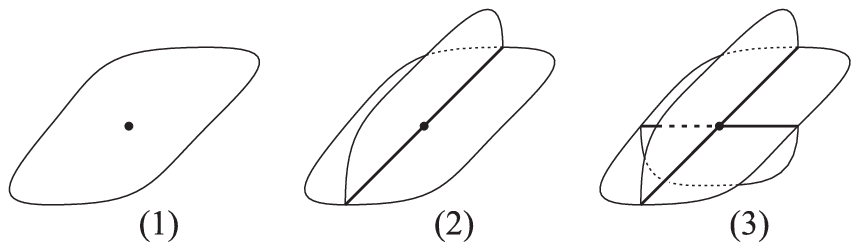}
\nota{Neighbourhoods of points in a standard polyhedron.} 
\label{standard_nhbds:fig}
\end{center}
\end{figure}

\paragraph{Definition of complexity}
A compact 2-dimensional
polyhedron $P$ is called \emph{simple} if the link of every point in $P$ is contained in
the 1-skeleton $\Delta^{(1)}$ of the tetrahedron. A point, a compact graph, a
compact surface are thus simple. Three important possible kinds of
neighbourhoods of points are shown in Fig.~\ref{standard_nhbds:fig}. 
A point having the whole of $\Delta^{(1)}$ as a link is 
called a \emph{vertex}, and its regular neighbourhood
is as shown in Fig.~\ref{standard_nhbds:fig}-(3). 
The set $V(P)$ of the vertices of $P$ consists of isolated points, so
it is finite. Points, graphs and surfaces of course do not contain vertices.
A compact polyhedron $P$ contained in the interior of a compact manifold $M$ 
with $\partial M\neq\emptyset$ is a \emph{spine} of $M$
if $M$ collapses onto $P$, \emph{i.e.}~if $M\setminus P \cong \partial M \times [0,1)$.
The \emph{complexity} $c(M)$ of
a 3-manifold $M$ is now defined as the minimal number of vertices of
a simple spine of either $M$, when $\partial M\ne\emptyset$, or $M$ minus some balls,
when $M$ is closed.

Since a point is a spine of the ball, a graph is a spine of a handlebody, and a surface
is a spine of an interval bundle, and
these spines do not contain vertices, the corresponding manifolds have
complexity zero. This shows that $c$ is not finite-to-one on manifolds
containing essential discs or annuli.

In general, to compute the
complexity of a manifold one must look for its \emph{minimal} spines,
\emph{i.e.}~the simple spines with the lowest number of vertices.
It turns out~\cite{Mat, Mat-Exp} that $M$ is $\partial$-irreducible and acylindrical if and only
if it has a minimal spine
which is \emph{standard}. A polyhedron is standard when
every point has a neighbourhood of one of the types (1)-(3) shown in
Fig.~\ref{standard_nhbds:fig}, and the sets of such points induce a
cellularization of $P$. That is, defining $S(P)$ as the set of points of type (2) or (3),
the components of $P \setminus S(P)$ 
should be open discs -- the \emph{faces} -- and the components of $S(P)\setminus V(P)$ 
should be open
segments -- the \emph{edges}. 

\begin{figure}
\begin{center}
\mettifig{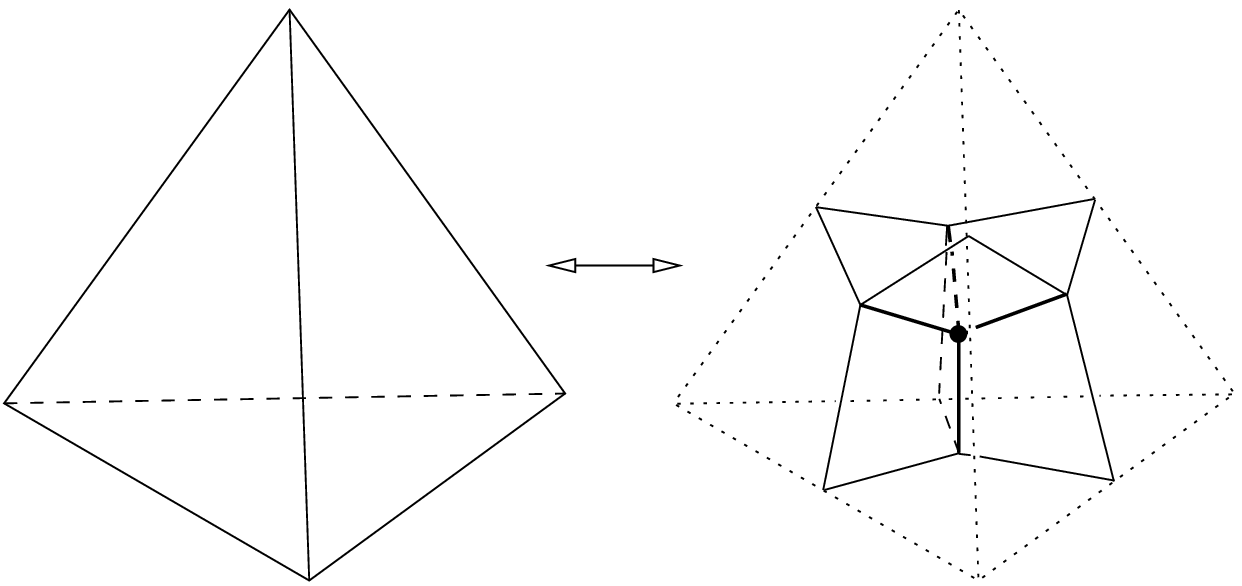, width=8 cm}
\nota{Duality between ideal triangulations and standard spines.} \label{dualspine:fig}
\end{center}
\end{figure}

The spines we are interested in are therefore standard and minimal. 
A standard spine is naturally dual to an ideal triangulation of $M$, 
as suggested in Fig.~\ref{dualspine:fig}. 
Moreover, by definition of $\widetilde\calH$ and the results of Matveev just cited, 
a manifold $M$ with $\chi(M)<0$ belongs to
$\widetilde\calH$ if and only if it has a standard minimal spine. 
These two facts imply the assertion
already stated that $c=t$ on $\widetilde\calH$ and $c<t$ 
outside $\widetilde\calH$ on manifolds with negative $\chi$.

\begin{figure}
\begin{center}
\mettifig{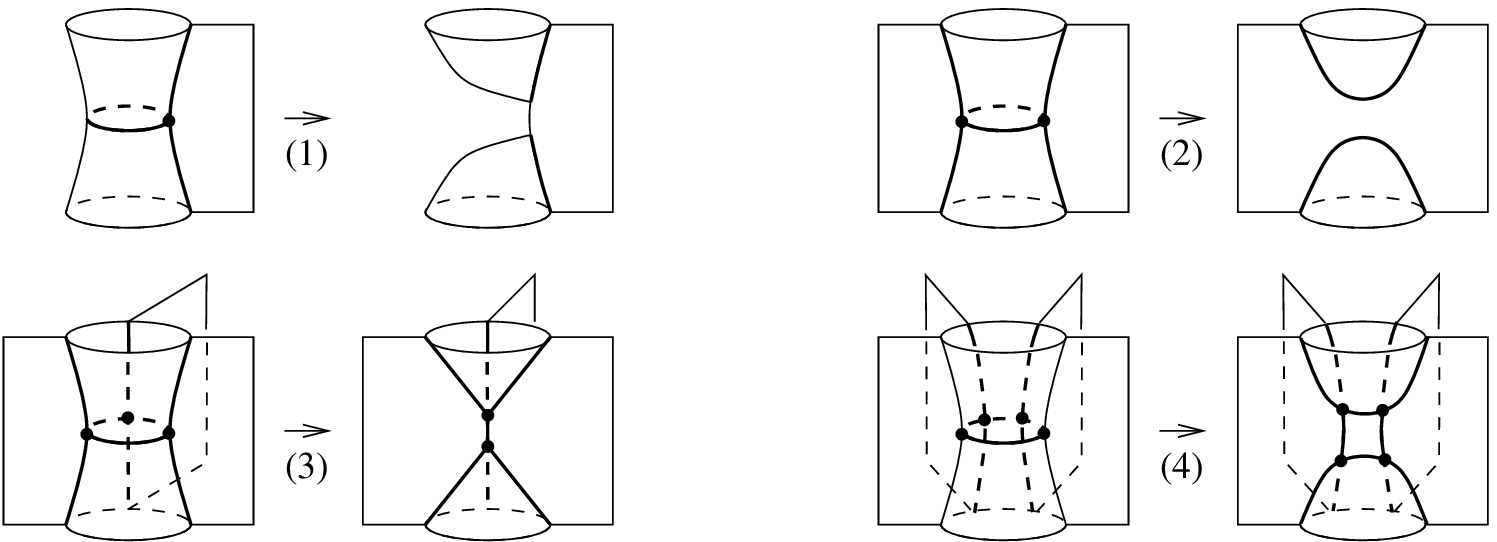, width=12 cm}
\nota{Moves on simple spines.} \label{moves:fig}
\end{center}
\end{figure}

\paragraph{Enumeration}
A naive approach to 
the classification of all manifolds in $\widetilde\calH_n$ for a fixed $n$
would be as follows:
\begin{enumerate}
\item Construct the finite
list of all standard polyhedra 
with $n$ vertices that are spines of some
orientable manifold (each such polyhedron is the spine of a unique manifold);
\item Check which of these spines are minimal, and discard the non-minimal ones;
\item Compare the corresponding manifolds for equality.
\end{enumerate}
Step (1) is feasible (even if the resulting list is very long), but step (2) is not, because
there is no general algorithm to tell if a given spine is minimal or not. In our classification
of $\widetilde\calH_3$ and $\widetilde\calH_4$ we have only performed some 
\emph{minimality tests}, and we have actually used them during the construction of the list,
to cut the ``dead branches'' at their bases and hence get not too huge a list. Our tests
are based on the moves shown in Fig.~\ref{moves:fig}, 
which are easily seen to transform a spine
of a manifold into another spine of the same manifold.  Namely, we have used the
following fact:
\begin{itemize}
\item If a spine $P$ of the list transforms into another one with less than $n$ vertices
via a combination of the moves of Fig.~\ref{moves:fig}, then $P$ is not minimal so it can 
be discarded.
\end{itemize}

\begin{rem}
\emph{Starting from a standard spine,
move (1) of Fig.~\ref{moves:fig} always leads to a simple but non-standard spine,
and move (2) also does on some spines, whereas moves (3) and (4) always give standard
spines. In particular, only moves (3) and (4) have counterparts at the level of
triangulations. This extra flexibility of simple spines compared to triangulations
is crucial for the enumeration.}
\end{rem}

Having obtained a list of candidate minimal spines with $n$ vertices, we conclude
the classification of $\widetilde\calH_n$ for $n=3,4$ as follows:
\begin{itemize}
\item For each spine in the list we write and try to solve numerically
the hyperbolicity equations, and if we find a solution we compute the canonical
Kojima decomposition, as discussed in Section~\ref{hyperbolicity:section}.
Solutions are found in all cases for $n=3$ and in all but $6$ cases for $n=4$.
All $6$ non-hyperbolic spines contain Klein bottles, so the corresponding 
manifolds cannot be hyperbolic;
\item Comparing the canonical decompositions of the hyperbolic manifolds thus found
and making sure they do not belong to $\calH_m$ for $m<n$,
we classify $\calH_n$. This gives $\widetilde\calH_3=\calH_3$ and $\calH_4$;
\item We show that the $6$ non-hyperbolic spines give distinct manifolds whose
complexity cannot be less than $4$, proving that $\widetilde\calH_4\setminus\calH_4$
contains $6$ manifolds.
\end{itemize}
The rest of this section is devoted to proving the last step and the 
assertions of Remark~\ref{toroidal:rem}.

\begin{figure}
\begin{center}
\mettifig{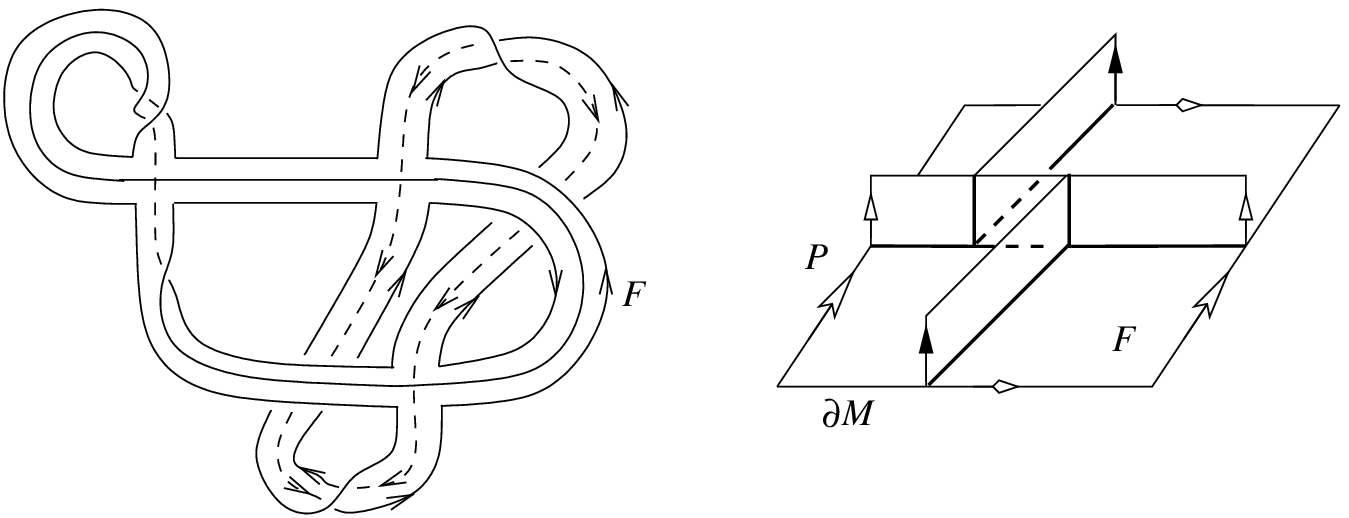, width=10 cm}
\nota{Left: a regular neighbourhood of $S(P)$; the rest of $P$ is obtained by
attaching two discs. Right: a regular neighbourhood in $P$ of the torus $T=\overline{F}$;
arrows indicate gluings.} \label{cusp:fig}
\end{center}
\end{figure}

\paragraph{Classification of $\widetilde\calH_4\setminus\calH_4$}
To analyze the $6$ non-hyperbolic spines with $4$ vertices we need more information
on the cusped element $M$ of $\calH_3$. Its unique minimal spine $P$ (described
in Fig.~\ref{cusp:fig}-left) has two faces, one of which, denoted by $F$,
is an open hexagon whose closure
in $P$ is a torus $T$. Since a neighbourhood of $T$ in $P$ is as 
in Fig.~\ref{cusp:fig}-right,
$P\setminus F$ is incident to $T$ on one side. Moreover the cusp of $M$ lies on the
other side of $T$, so $T$ can be viewed as the torus boundary component of the 
compactification of $M$.

\begin{figure}
\begin{center}
\mettifig{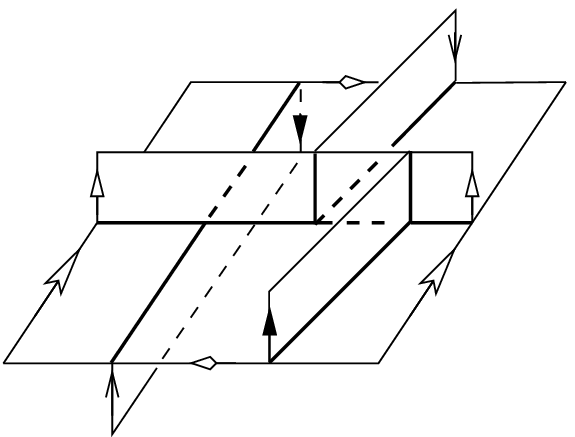, width=4 cm}
\nota{A simple polyhedron with $\theta$-shaped boundary.} \label{KtildeI:fig}
\end{center}
\end{figure}

Let us now consider the polyhedron $Q$ of Fig.~\ref{KtildeI:fig}, that one easily
sees to be a spine of the twisted 
interval bundle $K\timtil I$ over the Klein bottle. 
Note also that $Q$ has a natural $\theta$-shaped
boundary $\partial Q$ (a graph with two vertices and three edges)
that we can assume to lie on $\partial(K\timtil I)$. Now, if $P$ and $F$ are those
of Fig.~\ref{cusp:fig}, $P\setminus F$ also has a $\theta$-shaped boundary, and it turns
out that all the $6$ non-hyperbolic candidate minimal spines with $4$ vertices
have the form $(P\setminus F)\cup_\psi Q$, for some homeomorphism
$\psi:\partial Q\to\partial(P\setminus F)$. It easily follows that the 
associated manifold is $M\cup_\Psi(K\timtil I)$ where $\Psi:\partial(K\timtil I)\to T$
is the only homeomorphism extending $\psi$. 

Let us now choose a homology basis on 
$\partial(K\timtil I)$ so that the three slopes contained in
$\partial Q$ are $0,1,\infty\in\matQ\cup\{\infty\}$.
Doing the same on $T$ we see that $\Psi$ must map
$\{0,1,\infty\}$ to itself, so its matrix in ${\rm GL}_2(\matZ)$ must be
one of the following $12$ ones:
$$\pm\matr 1001,\qquad \pm\matr{-1}0{-1}1,\qquad \pm\matr 1{-1}0{-1},$$ 
$$\pm\matr 0110,\qquad \pm\matr{-1}1{-1}0,\qquad \pm\matr 0{-1}1{-1}.$$
Moreover the $6$ spines in question realize up to sign all these matrices.
Now the JSJ decomposition of $M\cup_\Psi(K\timtil I)$
consists of $M$ and $K\timtil I$, so $M\cup_\Psi(K\timtil I)$ is classified by the
equivalence class of $\Psi$ under the action of the automorphisms of $M$
and $K\timtil I$. But $M$ has no automorphisms, and 
it is easily seen that the only automorphism of
$K\timtil I$ acts as minus the identity on $\partial(K\timtil I)$. Therefore the $6$ spines
represent different manifolds. Moreover they are $\partial$-irreducible,
acylindrical, and non-hyperbolic, so they cannot belong to
$\widetilde\calH_m$ for $m<4$, and the classification is complete.

\section{Hyperbolicity equations and the tilt
formula}\label{hyperbolicity:section}
In this section we recall how an ideal triangulation can be used
to construct a hyperbolic structure with geodesic boundary
on a manifold, and how an ideal triangulation can be promoted to
become the canonical Kojima decomposition of the manifold.
We first treat the compact case and then sketch the variations
needed for the case where also some cusps exist. For all details
and proofs (and for some very natural terminology that we use 
here without giving actual definitions) we address the reader to~\cite{FriPe}.

\paragraph{Moduli and equations}
The basic idea for constructing a hyperbolic structure via an ideal
triangulation is to realize the tetrahedra as special geometric blocks
in $\matH^3$ and then require that the structures match when the
blocks are glued together. To describe the blocks to be used we first 
define a truncated tetrahedron $\Delta\!^*$ as a tetrahedron minus open stars of its
vertices. Then we call \emph{hyperbolic truncated tetrahedron} a realization of
$\Delta\!^*$ in $\matH^3$ such that the truncation triangles
and the lateral faces of $\Delta\!^*$ are geodesic triangles and
hexagons respectively, and the dihedral angle between a triangle
and a hexagon is always $\pi/2$.
Now one can show that:
\begin{itemize}
\item A hyperbolic structure on a combinatorial truncated 
tetrahedron is determined by the 6-tuple of dihedral angles along the
internal edges; 
\item The only restriction on this 6-tuple of positive
reals comes from the fact that the angles of each
of the four truncation triangles sum up to less than $\pi$;
\item The lengths of the internal edges can be computed as explicit 
functions of the dihedral angles;
\item A choice of hyperbolic structures on the tetrahedra of an ideal
triangulation of a manifold $M$ gives rise to a hyperbolic structure on $M$
if and only if all matching edges have the same length and the total dihedral
angle around each edge of $M$ is $2\pi$.
\end{itemize}
Given a triangulation of $M$ consisting of $n$ tetrahedra one then has
the \emph{hyperbolicity equations}: a system of $6n$ equations with unknown
varying in an open set of $\matR^{6n}$. We have solved these
equations using Newton's method with partial pivoting, after having explicitly
written the derivatives of the length function.

\paragraph{Canonical decomposition}
Epstein and Penner~\cite{epstein-penner}
have proved that cusped hyperbolic manifolds without 
boundary have a \emph{canonical decomposition}, and Kojima~\cite{kojima,Ko} 
has proved the same for
hyperbolic 3-manifolds with non-empty geodesic boundary.
This gives the following 
very powerful tool for recognizing manifolds: \emph{$M_1$ and $M_2$
are isometric (or, equivalently, homeomorphic) if and only
if their canonical decompositions are combinatorially equivalent}.
We have always checked equality and inequality of the manifolds 
in our census using this criterion, and we have proved that the cusped
element of $\calH_3$ has no non-trivial automorphism by showing 
that its canonical decomposition has no combinatorial automorphism.

Before explaining the lines along which we have found
the canonical decomposition of our manifolds, let us 
spend a few more words on the decomposition itself.
In the cusped case its blocks are 
\emph{ideal polyhedra}, whereas 
in the geodesic boundary case they
are \emph{hyperbolic truncated polyhedra} (an obvious generalization
of a truncated tetrahedron). In both 
cases the decomposition is obtained by projecting first to $\matH^3$ and then
to the manifold $M$ the faces of the convex hull of a certain family $\calP$
of points in Minkowsky 4-space. In the cusped case these points lie
on the light-cone, and they are the duals of the horoballs projecting in $M$
to Margulis neighbourhoods of the cusps. In the geodesic boundary case the
points lie on the hyperboloid of equation $\|x\|^2=+1$, and they
are the duals of the hyperplanes giving $\partial \widetilde{M}$, 
where $\widetilde{M}\subset\matH^3$ is a universal cover of $M$. 

\paragraph{Tilts}
Assume $M$ is a hyperbolic 3-manifold, either cusped without boundary or
compact with geodesic boundary, and let a geometric triangulation $\calT$ of 
$M$ be given. One natural issue is then
to decide if $\calT$ is the canonical decomposition
of $M$ and, if not, to promote $\calT$ to become canonical.
These matters are faced using the \emph{tilt formula}~\cite{weeks:tilt, Ush},
that we now describe.

If $\sigma$ is a $d$-simplex in $\calT$, the ends of its lifting 
to $\matH^3$ determine (depending on the nature of $M$) 
either $d+1$ Margulis horoballs 
or $d+1$ components of $\partial\widetilde{M}$, whence $d+1$ points
of $\calP$. Now let two tetrahedra $\Delta_1$ and $\Delta_2$ share a $2$-face
$F$, and let $\widetilde{\Delta}_1,\widetilde{\Delta}_2$ and $\widetilde F$ be 
liftings of $\Delta_1,\Delta_2$ 
and $F$ to $\matH^3$ such that $\widetilde{\Delta}_1\cap\widetilde{\Delta}_2=\widetilde{F}$.
Let $\overline{F}$ be
the $2$-subspace in Minkowsky 4-space that
contains the three points of $\calP$ determined by $\widetilde{F}$. For $i=1,2$ let 
$\overline{\Delta}^{(F)}_i$ be the half-$3$-subspace bounded by 
$\overline{F}$ and containing the fourth point of $\calP$ determined
by $\widetilde{\Delta}_i$. Then one can show that $\calT$ is canonical if and only if, 
whatever $F,\Delta_1,\Delta_2$, the convex hull of the
half-$3$-subspaces
$\overline{\Delta}^{(F)}_1$ and $\overline{\Delta}^{(F)}_2$
does not contain the origin of Minkowsky $4$-space, and
the half-$3$-subspaces themselves lie on 
distinct $3$-subspaces. Moreover, if
the first condition is met for all triples 
$F,\Delta_1,\Delta_2$, the canonical
decomposition is obtained by merging together the tetrahedra
along which the second condition is not met.

The tilt formula defines a real number $t(\Delta,F)$ describing the
``slope'' of $\overline{\Delta}^{(F)}$. More precisely,
one can translate the two conditions of the previous paragraph 
into the inequalities $t(\Delta_1,F)+t(\Delta_2,F)\leq 0$ and
$t(\Delta_1,F)+t(\Delta_2,F)\neq 0$ respectively.
Since we can compute tilts explicitly
in terms of dihedral angles, this gives a very efficient criterion
to determine whether $\calT$ is canonical or a subdivision
of the canonical decomposition. Even more, it suggests where to 
change $\calT$ in order to make it more likely to be canonical, 
namely along $2$-faces where the total tilt is positive.  This is achieved
by 2-to-3 moves along the offending faces, as 
discussed in~\cite{FriPe}. We only note here that the evolution
of a triangulation toward the canonical decomposition is not quite
sure to converge in general, but it always does in practice, and
it always did for us. We also mention that our computer program is only
able to handle triangulations: whenever some mixed negative and zero tilts
were found, the canonical decomposition was later worked out by hand.

\paragraph{Cusped manifolds with boundary}
When one is willing to accept both compact geodesic boundary and toric
cusps (but not annular cusps) the same strategy for constructing
the structure and finding the canonical decomposition applies, but many
subtleties and variations have to be taken into account. Let us quickly
mention which.

\emph{Moduli}. To parametrize tetrahedra one must consider that if a
vertex of some $\Delta$ lies in a cusp then the corresponding truncation triangle
actually disappears into an ideal vertex (a point of $\partial\matH^3$).
At the level of moduli this translates into the condition that the
triangle be Euclidean, \emph{i.e.} that its 
angles sum up to precisely $\pi$.

\emph{Equations}. If an internal edge ends in a cusp then its length is infinity,
so some of the length equations must be dismissed when there are cusps. On the other
hand, when an edge is infinite at both ends, one must make sure that the gluings
around the edge do not induce a sliding along the edge. This translates into
the condition that the \emph{similarity moduli} of the Euclidean triangles
around the edge have product $1$. This ensures consistency of the hyperbolic
structure, but one still has to impose completeness of cusps. Just as in the case
where there are cusps only, this amounts to requiring that the similarity tori
on the boundary be Euclidean, which translates into the \emph{holonomy equations}
involving the similarity moduli.

\emph{Canonical decomposition}. When there are cusps, 
the set of points $\calP$ to take the convex hull of consists of the
duals of the planes in $\partial\widetilde{M}$ and of some points on the
light-cone dual to the cusps. The precise discussion on how to 
choose these extra points is too complicated to be reproduced 
here (see~\cite{FriPe}), but the implementation of the choice was 
actually very easy in the (not many) cusped members of our census. 

The computation of tilts and the discussion on how to find 
the canonical decomposition are basically unaffected by the 
presence of cusps.

\newpage

\section*{Appendix: Tables of volumes}

\begin{table}[h]\begin{center}

\nota{Number of manifolds per value of volume for the compact elements of $\calH_3$ with 
boundary of genus $2$ and canonical decomposition into thee tetrahedra.}\label{c3:table}
\end{center}\end{table}

\begin{table}[hb]\begin{center}
%
\nota{Number of manifolds per value of volume for compact elements of $\calH_4$ with boundary $\Sigma_3$
and canonical decomposition into four tetrahedra -- first part.}\label{c4:bd3:K4:first:table}
\end{center}\end{table}

\begin{table}\begin{center}
%
\nota{Number of manifolds per value of volume for compact elements of $\calH_4$ with boundary $\Sigma_3$
and canonical decomposition into four tetrahedra -- second part. Note the changes of scale.}\label{c4:bd3:K4:second:table}
\end{center}\end{table}

\begin{table}\begin{center}
%
\nota{Number of manifolds per value of volume for compact elements of $\calH_4$ with boundary $\Sigma_2$ 
and canonical decomposition into four tetrahedra -- first part.}\label{c4:bd2:K4:first:table}
\end{center}\end{table}

\begin{table}\begin{center}
%
\nota{Number of manifolds per value of volume for compact elements of $\calH_4$ with boundary $\Sigma_2$ 
and canonical decomposition into four tetrahedra -- second part.
Note the change of scale on volumes.}\label{c4:bd2:K4:second:table}
\end{center}\end{table}

\begin{table}\begin{center}
%
\nota{Number of manifolds per value of volume for compact elements of $\calH_4$ with boundary $\Sigma_2$ 
and canonical decomposition into five tetrahedra.}\label{c4:bd2:K5:table}
\end{center}\end{table}

\begin{table}\begin{center}
%
\nota{Number of manifolds per value of volume for one-cusped elements of $\calH_4$
with boundary $\Sigma_2$
and canonical decomposition into four tetrahedra.}\label{c4:cusp:table}
\end{center}\end{table}

\vspace{.7 cm}

\noindent Scuola Normale Superiore, Piazza dei Cavalieri 7, 56127 Pisa, Italy\\ 
frigerio@sns.it
\vspace{.2 cm}

\noindent
Dipartimento di Matematica, 
Via F. Buonarroti 2, 
56127 Pisa, Italy\\ 
martelli@mail.dm.unipi.it
\vspace{.2 cm} 

\noindent
Dipartimento di Matematica Applicata,
Via Bonanno Pisano 25B,
56126 Pisa, Italy\\ 
petronio@dm.unipi.it

\end{document}